\documentclass[a4paper,12pt]{article}

\usepackage{amsfonts, amsmath, amssymb}

\newtheorem{theorem}{Theorem}
\newtheorem{lemma}[theorem]{Lemma}
\newtheorem{proposition}[theorem]{Proposition}
\numberwithin{theorem}{section}

\begin{document}

\title{\bf On special values of standard $L$-functions of Siegel cusp eigenforms of genus 3}



\author{Anh Tuan DO
 \footnote{Institut Fourier, Universit\'{e} Grenoble I, 38402 Saint Martin D'H\`{e}res, France, {\tt anh-tuan.do@ujf-grenoble.fr}}
 \and Kirill VANKOV\footnote{{\tt kirill.vankov@gmail.com}}
}

\date{December 20, 2013}

\maketitle

\begin{abstract}
 We explicitly compute the special values of the standard $L$-function $L(s, F_{12}, \mathrm{St})$ at 
 the critical points $s\in\{-8, -6, -4, -2, 0, 1, 3, 5, 7, 9\}$, where $F_{12}$ is the unique (up to a scalar) 
 Siegel cusp form of degree $3$ and weight $12$, which was constructed by Miyawaki. 
 These values are proportional to the product of the Petersson norms of symmetric square of Ramanujan's 
 $\Delta$ and the cusp form of weight $20$ for ${\rm SL}_2(\mathbb{Z})$ by a rational number and some power of $\pi$.  
 We use the Rankin-Selberg method and apply the Holomorphic projection to compute these values. 
 To our knowledge this is the first example of a standard $L$-function of Siegel cusp form of degree $3$, 
 when the special values can be computed explicitly.
\end{abstract}

\section{ Introduction}
In the remarkable paper \cite{miyawaki:1992} I.~Miyawaki considered certain Siegel cusp forms of degree $3$, 
and on the basis of some numerical calculations, he was able to make interesting conjectures about the degeneration of the standard and spinor $L$-functions associated to such cusp forms. Several years later T.~Ikeda \cite{ikeda:2006} proved Miyawaki's conjecture related to the standard $L$-function. Basically, he was able to construct an explicit lifting from Siegel cusp forms of degree $r$ to Siegel cusp forms of degree $r+2n$. In particular, it turns out that the cusp form of degree $3$ and weight $12$ is a basic example of this lifting (for $r=1, n=1$).

Recall that Miyawaki constructed his numerical examples by means of theta functions with spherical functions.  
Namely, let $E_8$ be the unique even unimodular lattice of rank $8$ i.e.,
\begin{equation}
 E_8=\left\{ {}^t(x_1,\dotsc,x_8)\in\mathbb{R}^8 \left|\begin{aligned}&{2x_i\in\mathbb{Z} (i=1,\dotsc,8),}\\&{x_1+\dotsb+x_8\in2\mathbb{Z},}\\&{x_i-x_j\in\mathbb{Z}}\end{aligned}\right. \right\}\,,
\end{equation}
and
\begin{equation}
Q=\begin{pmatrix}1&0&0&i&0&0&0&0\\
                 0&1&0&0&i&0&0&0\\
                 0&0&1&0&0&i&0&0
       \end{pmatrix}
\end{equation}
be $3\times8$ matrix.  Then the theta series
\begin{align}
F_{12}(Z)&=\sum_{v_1,v_2,v_3\in E_8}\Re\left(\det(Q\cdot(v_1,v_2,v_3))^8\right)\,\exp\left(\pi i\sigma(\langle v_i,v_j\rangle)Z)\right)\\
&=\sum_{N>0}a(N)\exp(2\pi i \sigma(NZ))
\end{align}
is a cusp form of weight $12$ with respect to $\mathrm{Sp}_3(\mathbb{Z})$, where $(\left<v_i,v_j\right>)$ 
denotes the matrix composed by the entries, which are the scalar product of two vectors $v_i$ and $v_j$, 
$Z\in\mathfrak{H}^3$, $\sigma$ is the trace operator, and
\begin{align}
&\mathrm{Sp}_n(\mathbb{Z})=\left\{ M\in\mathrm{Mat}_{2n}(\mathbb{Z})\,\mid\,M\,\left(\begin{smallmatrix}0&\mathrm{I}_n\\-\mathrm{I}_n&0\end{smallmatrix}\right)^t\!M=\left(\begin{smallmatrix}0&\mathrm{I}_n\\-\mathrm{I}_n&0\end{smallmatrix}\right) \right\}\,,\\
&\mathfrak{H}^3=\left\{Z={}^t\!Z=X+iY\,\mid\,X,Y\in \mathrm{Mat}_3(\mathbb{R}),\,Y>0\right\}\,.
\end{align}

Therefore, the explicit Fourier expansion can be computed
\begin{align}
F_{12}&=1\,\cdot\,q^{\left(\begin{smallmatrix}1&\frac{1}{2}&\frac{1}{2}\\\frac{1}{2}&1&\frac{1}{2}\\\frac{1}{2}&\frac{1}{2}&1\end{smallmatrix}\right)}
+164\,\cdot\,q^{\left(\begin{smallmatrix}1&0&0\\0&1&0\\0&0&1\end{smallmatrix}\right)}
+1328\,\cdot\,q^{\left(\begin{smallmatrix}3&0&\frac{1}{2}\\0&1&\frac{1}{2}\\\frac{1}{2}&\frac{1}{2}&1\end{smallmatrix}\right)}
\\
&-1008\,\cdot\,q^{\left(\begin{smallmatrix}2&0&0\\0&1&0\\0&0&1\end{smallmatrix}\right)}
-131776\,\cdot\,q^{\left(\begin{smallmatrix}2&1&1\\1&2&1\\1&1&2\end{smallmatrix}\right)}
-6816512\,\cdot\,q^{\left(\begin{smallmatrix}2&0&0\\0&2&0\\0&0&2\end{smallmatrix}\right)}
+\dotsc
\end{align}

The purpose of this note is to show that at each critical point $s$ one can provide the explicit rational number $R(s)$ 
and power of $\pi$ such that
\begin{align}
\label{eq1.2}
L(s, F_{12}, \mathrm{St})=R(s)\,\pi^{\alpha_s}\,\langle \Delta, \Delta\rangle\,\langle g_{20}, g_{20}\rangle,
\end{align}
where $\Delta$ is Ramanujan's discriminant cusp form and $g_{20}$ is the cusp form of weight $20$ of level $1$.

The functional equation for a standard $L$ function in general case is proved by B\"{o}cherer \cite{bocherer:1985}.  
For $F_{12}$ it is as follows:
\begin{equation}
\Psi(s,F_{12},\mathrm{St})=\Psi(1-s,F_{12},\mathrm{St})\,,
\end{equation}
where
\begin{align}
&\Psi(s,F_{12},\mathrm{St})=\gamma(s)\;L(s,F_{12},\mathrm{St})\,,\\
&\gamma(s)=C\;2^{-3s}\;\pi^{-7s/2}\;\Gamma\left(\textstyle\frac{s+1}{2}\right)\;\Gamma(s+9)\;\Gamma(s+10)\;\Gamma(s+11)
\end{align}
with some non-zero constant $C$. Critical points for $L(s,F_{12},\mathrm{St})$ in the sense of Deligne \cite{deligne:1979} 
are $s\in\{-8,-6,-4,-2,0,1,3,5,7,9\}$. In particular, gamma factor $\gamma(s)$ on the left hand side of the functional 
equation and $\gamma(1-s)$ on the right hand side is finite at these points.

Due to Miyawaki and Ikeda
\begin{equation}
\label{eq1.1}
L(s,F_{12},\mathrm{St})=L(s+11,\Delta\otimes\Delta)\,L(s+10,g_{20})\,L(s+9,g_{20})\,.
\end{equation}
This theorem is the starting point of our investigation.  We compute the result in two steps: first, 
compute $L(s+11,\Delta\otimes\Delta)$ for all points we are interested in, and then compute $L(s+10,g_{20})\,L(s+9,g_{20})$. 
We use Rankin-Selberg method to represent the $L$-funcion of multiplicative convolution as an integral over a fundamental domain, 
which can be written in terms of the Petersson scalar product. And, in order to compute it, we use a holomorphic projection.

To our knowledge this is the first example of a standard $L$-function of Siegel cusp form of degree $3$, 
when the special values can be computed explicitly. The case of spinor $L$-function was treated in our other 
work \cite{chiera-vankov:2008}. 
The authors are very grateful to Alexei Panchishkin for intensive and encouraging discussions, 
this work was undertaken due to his keen interest into the subject. 

\section{Generalities and notation} 

Let $\mathfrak{H}=\{z\in\mathbb{C}\mid\Im(z)>0\}$ be the upper half-plane. For a positive integer $k$ and a Dirichlet 
character $\chi$ modulo a positive integer $N$ such that $\chi(-1)=(-1)^k$, we denote by $\mathcal{M}_k(\Gamma_0(N),\chi)$ 
the vector space of all holomorphic modular forms $f(z)$ of weight $k$ satisfying
\begin{align}
f(\gamma(z))=\chi(d)(cz+d)^kf(z)\mbox{ for all } \gamma=\begin{pmatrix}a&b\\c&d\end{pmatrix}\in\Gamma_0(N),
\end{align}
where the variable $z\in\mathfrak{H}$, $\gamma(z)=\displaystyle\frac{az+b}{cz+d}$, and
\begin{align}
\Gamma_0(N)=\left\{\begin{pmatrix}a&b\\c&d\end{pmatrix}\in{\rm SL}_2(\mathbb{Z})\mid c\equiv 0\mod N\right\}.
\end{align}
We denote by $S_k(N,\chi)$ the subspace of $\mathcal{M}_k(\Gamma_0(N),\chi)$ consisting of all cusp forms.  
Every element $f$ of $\mathcal{M}_k(\Gamma_0(N),\chi)$ has a Fourier expansion
\begin{align}
f(z)=\sum_{n=0}^\infty a(n)\,q^n,
\end{align}
where $q=\exp(2\pi iz)$ and $a(n)$ are complex numbers in general.

The $L$-function associated to $f$ is defined as
\begin{equation}
L(s,f)=\sum_{n=1}^\infty a(n)\,n^{-s}. 
\end{equation}
More generally, with an arbitrary Dirichlet character $\omega$, the twisted $L$-function is defined as 
$\displaystyle L(s,f,\omega)=\sum_{n=1}^\infty a(n)\,\omega(n)\,n^{-s}$. These $L$-functions can be also 
written in the form of Euler product:
\begin{align}
L(s,f)&=\prod_{p\text{ prime}}(1-a(p)\,p^{-s}+\chi(p)\,p^{k-1-2s})^{-1}\,,\\
L(s,f,\omega)&=\prod_{p\text{ prime}}(1-a(p)\,\omega(p)\,p^{-s}+\chi(p)\,\omega(p)^2\,p^{k-1-2s})^{-1}\,.
\end{align}
Next, the Dirichlet $L$-series for any character $\chi$ of conductor $N$ is given by
\begin{equation}
\label{eq:Dirichlet_L-series}
L(s,\chi)=\sum_{n=1}^\infty\chi(n)\,n^{-s}
\end{equation}
and its Euler product
\begin{equation}
\label{eq:Euler_product_Dirichlet_L-series}
L(s,\chi)=\prod_{p\,\nmid N}\frac{1}{1-\chi(p)\,p^{-s}}\,.
\end{equation}
In the case, when $\chi$ is the identity Dirichlet character, the latter series is Riemann's zeta function 
$\displaystyle\zeta(s)=\sum_{n=1}^\infty\frac{1}{n^s}$.

Let $\displaystyle g(z)=\sum_{n=0}^\infty b(n)\,q^n \in \mathcal{M}_l(\Gamma_0(N),\xi)$ be another modular 
form of weight $l$ with Fourier coefficients $b(n)$. The $L$ function associated to two modular forms $f$ 
and $g$ is given by the additive convolution
\begin{equation}
L(s,f,g)=\sum_{n=1}^\infty a(n)\,b(n)\,n^{-s}\,.
\end{equation}

Another type of $L$-functions associated to two modular forms is Rankin's product $L$-function (multiplicative convolution).  
It is denoted by \mbox{$L(s,f\otimes g)$} and defined as (see \cite[page 786]{shimura:1976}):
\begin{align}
\label{eq2.1}
L(s, f\otimes g)=L_N(2s+2-k-l, \chi\xi)L(s, f, g),
\end{align}
where $L$-function $L_N(s,\omega)$ with a Dirichlet character $\omega$ modulo $N$ is defined, as usual, 
in (\ref{eq:Dirichlet_L-series}) with $\omega(n)=0$ for $(n,N)\neq 1$, and the Euler factors in 
(\ref{eq:Euler_product_Dirichlet_L-series}), corresponding to the prime divisors of a number $N$, have been omitted.  
Note, that in the case, when $f$ and $g$ are cusp eigenforms, the left-hand side of (\ref{eq2.1}) is an Euler product 
of degree $4$ in view of the following Lemma:
\begin{lemma}[Lemma 1, \cite{shimura:1976}]
\label{lemma:Rankin}
Suppose we have formally
\begin{align}
\sum_{n=1}^\infty A(n)\,n^{-s}=\prod_p\left[(1-\alpha_p\,p^{-s})(1-\alpha_p^\prime\,p^{-s})\right]^{-1}
\end{align}
and
\begin{align}
\sum_{n=1}^\infty B(n)\,n^{-s}=\prod_p\left[(1-\beta_p\,p^{-s})(1-\beta_p^\prime\,p^{-s})\right]^{-1}\,.
\end{align}
Then
\begin{align}
&\sum_{n=1}^\infty A(n)\,B(n)\,n^{-s}=\\
&\quad\prod_p\frac{1-\alpha\alpha^\prime\beta\beta^\prime p^{-2s}}{(1-\alpha\beta p^{-s})(1-\alpha\beta^\prime p^{-s})(1-\alpha^\prime\beta p^{-s})(1-\alpha^\prime\beta^\prime p^{-s})}\,.
\end{align}
\end{lemma}

In the case when $f=g$, we introduce the convoluted zeta-function of $f$ with itself by defining Euler $p$-factor as
\begin{align}
L_p(s,f\otimes f)&=(1-\alpha_1^2p^{-s})^{-1}\,(1-\alpha_1\alpha_2p^{-s})^{-1}\,(1-\alpha_2^2p^{-s})^{-1}\\
&=(1-p^{-s+k-1})^{-1}\,(1-\alpha_1^2p^{-s})^{-1}\,(1-\alpha_2^2p^{-s})^{-1}\,.
\end{align}
Then the zeta-function $L(s,f\otimes f)$ is defined by
\begin{equation}
L(s,f\otimes f)=\prod_pL_p(s,f\otimes f)\,.
\end{equation}
Therefore, it is related to the additive convolution by
\begin{equation}
\label{eq:L(s,fxf)}
L(s,f\otimes f)=\frac{\zeta(2s+2-2k)}{\zeta(s+1-k)}\sum_{n=1}^\infty a^2(n)\,n^{-s}\,.
\end{equation}

For two elements $f,h\in\mathcal{M}_k(\Gamma_0(N))$ such that the product $fh$ is a cusp form, the Petersson inner 
product $\langle f, h\rangle$ is defined as
\begin{equation}
\label{eq:Petersson}
\langle f, h\rangle=\dfrac{1}{[\mathrm{SL}_2(\mathbb{Z}):\Gamma_0(N)]}\int_{\Phi_N}\overline{f(z)}\,h(z)\,y^{k-2}\,dx\,dy\,,
\end{equation}
where $z=x+iy$, $\Phi_N$ is a fundamental domain for $\mathfrak{H}$ modulo $\Gamma_0(N)$ and the bar denotes the complex conjugate. 
We also define $\left<f,h\right>$ by (\ref{eq:Petersson}) for nearly holomorphic modular forms $f$ and $h$ on $\mathfrak{H}$ 
whenever the integral is convergent (see \cite[section 8.2]{shimura:2007} for definition and properties of the nearly holomorphic 
modular forms).

We also briefly recall the definition of standard $L$-function for modular forms of arbitrary genus. 
Let $f\in\mathcal{M}_k^n(N, \psi)$ be an eigenfunction of all Hecke operators $f\longmapsto f|T, T\in\mathcal{L}_q^n(N)$ 
with $q$ being prime number, $q\nmid N$, so that $f|T=\lambda_f(T)f$. Then the number $\lambda_f(T)\in\mathbb{C}$ 
define a homomorphism $\lambda_f: \mathcal{L}\longrightarrow \mathbb{C}$ which is uniquely determined by $(n+1)$-tuple of numbers
\begin{align*}
(\alpha_0, \alpha_1, \cdots, \alpha_n)=(\alpha_{0, f}(q), \alpha_{1, f}(q), \cdots, \alpha_{n, f}(q))\in[(\mathbb{C})^{n+1}]^{W_n}
\end{align*}
which are called the Satake $q$-parameters of the modular form $f$.

Now let the variables $T_0, T_1, \cdots, T_n$ be equal to the corresponding Satake $q$-parameters 
$\alpha_{0, f}(q), \alpha_{1, f}(q), \cdots, \alpha_{n, f}(q)$ then
\begin{align*}
R_{f, q}(z)=\prod_{i=1}^n(1-\alpha_i^{-1}z)(1-\alpha_iz)\in\mathbb{Q}[\alpha_0^{\pm 1}, \cdots, \alpha_n^{\pm 1}].
\end{align*}
The standard zeta function of $f$ is defined by means of the Satake $p$-parameters as the following Euler product:
\begin{align}
L(s, f, \mathrm{St}, \chi)=\prod_{q\nmid N}L^{(q)}(s, f, \mathrm{St}, \chi) 
\end{align}
with 
\begin{align}
&L^{(q)}(s, f, \mathrm{St}, \chi)=(1-\chi(q)\psi(q)q^{-s})^{-1}R_{f, q}(\chi(q)\psi(q)q^{-s})^{-1}, \\
&L(s, f, \mathrm{St}, \chi)=\nonumber\\
&\qquad\prod_p\left(\left(1-\frac{\chi(p)}{p^s}\right)\prod_{i=1}^m\left(1-\frac{\chi(p)\alpha_i(p)}{p^s}\right)\left(1-\frac{\chi(p)\alpha_i(p)^{-1}}{p^s}\right)\right)^{-1}.
\end{align}

\section{Computation of $L(s,\Delta\otimes\Delta)$}
\label{sec:Delta}

The first term $L(s+11,\Delta\otimes\Delta)$ is the symmetric square of cusp form $\Delta$ by the mean of (\ref{eq2.1}). 
This $L$ function has been already studied in details by Rankin, Zagier, Li, Sturm and others. 
Using formulas from \cite{rankin:1939} and \cite[page 116]{zagier:1977} we have the values of $L(s, \Delta\otimes\Delta)$ 
at the points $s\in\{12, 14, 16, 18, 20\}$. The values at points $s\in\{3, 5, 7, 9, 11\}$ can be obtained from 
the functional equation (see \cite{li:1979})
\begin{align}
D^{*}(s, \Delta)=D^{*}(23-s, \Delta),
\end{align}
where
\begin{align}
D^*(s,\Delta)&=2^{-s}\,\pi^{-3s/2}\,\Gamma(s)\,\Gamma\left(\textstyle\frac{s-10}{2}\right)\,L(s,\Delta\otimes\Delta).
\end{align}
The results are presented in Appendix \ref{app:L(s,Delta2)}.

\section{The expression for $L(s, g_{20})L(s-1, g_{20})$}
\label{sec:L(s,g20)L(s-1,g20)}

We compute the product $L(s,g_{20})\,L(s-1,g_{20})$ at points critical points of $L(s, F_{12}, \mathrm{St})$ 
(note the shift by 10) $s\in\{2,4,6,8,10,11,13,15,17,19\}$ in two steps. The idea is to present the product 
$L(s,f)\,L(s-1,f)$ as an $L$ function of the Rankin convolution of $f$ with appropriate Eisenstein series and 
to use Rankin-Selberg method to relate the obtained expression to the Petersson inner product.

Let $g_{20}\in S_{20}$ be the cusp form of weight $20$ of level $1$:
\begin{align}
g_{20}(z)&=\sum_{n=1}^{\infty}a(n)q^n\\
&=q+456q^2+50652q^3-316352q^4+2377410q^5-\cdots
\end{align}
and its associated $L$-function is
\begin{align}
L(s, g_{20})=\sum_{n=1}^{\infty}a(n)n^{-s}.
\end{align}
Assume $1-a(p)X+p^{19}X^2=(1-\alpha_pX)(1-\alpha^\prime_pX)$, then $\alpha_p+\alpha^\prime_p=a(p)$, 
$\alpha_p\alpha^\prime_p=p^{19}$, and
\begin{align}
\label{eq:L(s,f)}
L(s, g_{20})=\prod_{p}((1-\alpha_pp^{-s})(1-\alpha^\prime_pp^{-s}))^{-1}.
\end{align}

Let
\begin{align}
G_2(z)&=-\frac{1}{24}+\sum_{n=1}^{\infty}\sigma_1(n)q^n\\
 & =-\frac{1}{24}+q+3q^2+4q^3+7q^4+\cdots,
\end{align}
where $\sigma_1(n)=\sum_{d|n}d$ is the divisor function defined as the sum of the divisors of $n$. 

Consider the Eisenstein series (see \cite[Corollary 7.2.14]{miyake:2006})
\begin{align}
G_{2, p}(z)&=G_2(z)-pG_2(pz)=\frac{p-1}{24}+\sum_{n=1}^{\infty}\sum_{d|n, p\nmid d}dq^n,
\end{align}
of weight $2$ for $\Gamma_0(p)$ and the corresponding $L$ series
\begin{align}
L(s,G_{2, p})&=\sum_{n=1}^\infty\sum_{\substack{d|n\\p\,\nmid\,d}}d\,n^{-s}\\
&=\sum_{\substack{d,d_1\geqslant 1\\p\,\nmid\,d}}d\,(d\,d_1)^{-s}\\
&=\sum_{\substack{d\geqslant 1\\p\,\nmid\,d}}d^{1-s}\sum_{d_1\geqslant 1}d_1^{-s}\\
&=(1-p^{1-s})\,\zeta(s-1)\,\zeta(s)\,.\\
\end{align}

Let us put $p=2$, then consider the following holomorphic modular form of weight $2$ and level $2$
\begin{align}
G_{2, 2}(z)&=G_2(z)-2G_2(2z)\\
&=\frac{1}{24}+q+q^2+4q^3+q^4+6q^5+4q^6+\cdots,
\end{align}
and
\begin{align}
L(s, G_{2, 2})=(1-2^{1-s})\,\zeta(s-1)\,\zeta(s).
\end{align}
Note, 
\begin{align}
G_{2, 2}=\sum_{n=0}^{\infty}b(n)q^n\in\mathcal{M}_2(\Gamma_0(2), \xi)
\end{align}
with
\begin{equation}
\xi(n)=(1~\mathrm{mod}~2)(n)=\begin{cases}1,&\text{if $n$ odd;}\\0,&\text{if $n$ even.}\end{cases}
\end{equation}

Similarly to \eqref{eq:L(s,f)}, consider Fourier coefficients and the decomposition
\begin{equation}
\begin{split}
L(s, G_{2, 2})&=\sum_{n=1}^{\infty}b(n)n^{-s}\\
&=(1-2^{1-s})\,\zeta(s-1)\,\zeta(s)\\
&=(1-2^{1-s})\prod_p((1-p^{1-s})(1-p^{-s}))^{-1}\\
&=\prod_p((1-\beta_pp^{-s})(1-\beta^\prime_pp^{-s}))^{-1},
\end{split}
\end{equation}
where $\beta(p)=1$ for all $p$, $\beta^\prime(2)=0$ and $\beta^\prime(p)=p$ for all odd primes.

Recall Lemma \ref{lemma:Rankin} and consider $g_{20}$ and $G_{2, 2}$ in its context. 
By definition (\ref{eq2.1}) and using Lemma \ref{lemma:Rankin} with $k=20$, $l=2$, $N=2$, $\chi=1$ we obtain
\begin{equation}
\begin{split}
L(&s, g_{20}\otimes G_{2, 2})=L_2(2s+2-20-2, \psi)L(s, g_{20}, G_{2, 2})\\
&=\prod_{p\not=2}(1-p^{20-2s})^{-1}\cdot\sum_{n=1}^{\infty}a(n)b(n)n^{-s}\\
&=\prod_{p\not=2}(1-p^{20-2s})^{-1}\\
&\quad\times\prod_{p}\frac{1-\alpha_p\alpha^\prime_p\beta_p\beta^\prime_pp^{-2s}}{(1-\alpha_p\beta_pp^{-s})(1-\alpha^\prime_p\beta_pp^{-s})(1-\alpha_p\beta^\prime_pp^{-s})(1-\alpha^\prime_p\beta^\prime_pp^{-s})}\\
&=\frac{1}{(1-\alpha_22^{-s})(1-\alpha^\prime_22^{-s})}\\
&\quad\times\prod_{p\not=2}\frac{1}{(1-\alpha_pp^{-s})(1-\alpha^\prime_pp^{-s})(1-\alpha_pp^{1-s})(1-\alpha^\prime_pp^{1-s})}\\
&=\prod_p\frac{1}{(1-\alpha_pp^{-s})(1-\alpha^\prime_pp^{-s})}\prod_{p\not=2}\frac{1}{(1-\alpha_pp^{1-s})(1-\alpha^\prime_pp^{1-s})}\\
&=L(s, g_{20})L(s-1, g_{20})(1-a(2)2^{1-s}+2^{19}2^{2-2s})\\
&=(1-456\cdot 2^{1-s}+2^{21-2s})L(s, g_{20})L(s-1, g_{20}).
\end{split}
\end{equation}

Finally, we obtain the following identity
\begin{align}\label{eq3.1}
L(s, g_{20})L(s-1, g_{20})=\frac{L(s, g_{20}\otimes G_{2, 2})}{(1-456\cdot 2^{1-s}+2^{21-2s})}.
\end{align}

\section{Computation of $L(s, g_{20}\otimes G_{2, 2})$}

Now we express $L(s, g_{20}\otimes G_{2, 2})$ (at integer points) as a multiple of Petersson inner product 
$\langle g_{20}, g_{20}\rangle$. Put $g_{20}^{*}(z)=\overline{g_{20}(-\bar{z})}=\sum \overline{a(n)}q^n$. 
We use Shimura's formula \cite[(2.4)]{shimura:1976}
\begin{align}\label{eq4.1}
L(s, &g_{20}\otimes G_{2, 2})=\\
&\frac{(4\pi )^s}{2\Gamma (s)}\int_{\Phi_2}\overline{g_{20}^{*}(z)}G_{(2, 2)}(z)E_{18, 2}(z, s-19, \xi ) y^{s-1}dxdy=\\
&\frac{(4\pi)^s}{2\Gamma(s)}\int_{\Phi_2}\overline{g_{20}^{*}(z)}G_{(2, 2)}(z)E_{18, 2}(z, s-19, \xi)y^{s-19}y^{18}dxdy=\\
&\frac{(4\pi)^s[{\rm SL}_2(\mathbb{Z}): \Gamma_0(2)]}{2\Gamma(s)}\langle g_{20}(z), G_{2, 2}(z)y^{19}E_{18, 2}(z, s-19, \xi)\rangle=\\
&\frac{3}{2}\frac{(4\pi)^s}{\Gamma(s)}\langle g_{20}(z), \mathcal{H}ol (G_{2, 2}(z)y^{19}E_{18, 2}(z, s-19, \xi))\rangle
 \end{align}
where $\langle f, g\rangle$ is the Petersson inner product (\ref{eq:Petersson}), $\Phi_2$ denotes a fundamental domain 
for $\Gamma_0(2)\backslash\mathfrak{H}$, $z=x+iy$,
\begin{equation}
E_{\lambda,N}(z,s,\xi)=\mathop{{\sum}^\prime}_{(m,n)} \xi(n)\,(mNz+n)^{-\lambda}|mNz+n|^{-2s}\,,
\end{equation}
$\mathop{{\sum}^\prime}$ denotes the summation over all $(m,n)\in\mathbb{Z}^2$, $(m,n)\neq(0,0)$, $\mathcal{H}ol(F)$ 
is the operator of holomorphic projection  (see \cite[(2.148)]{courtieu-panchishkin:2004}). It is defined so, 
that $\left<f,F\right>=\left<f,\mathcal{H}ol(F)\right>$ for all $f\in S_k(N,\psi)$.

In order to compute the holomorphic projection of the product in the last identity of (\ref{eq4.1}) 
$\mathcal{H}ol\left(G_{2,2}(z)\,(4\pi y)^{s+1-k}\,E_{k-2,2}(z,s+1-k,\xi)\right)$ we need the explicit Fourier 
coefficients of the product $G_{2,2}\,E_{k-2,2}$. The holomorphic projection belongs to the space of cups forms 
of weight $20$ and level $2$. It is a $4$-dimensional space, therefore we need at least four Fourier coefficients 
of it in order to identify it in a fixed basis. The Fourier coefficient of $G_{2,2}(z)$ are known. 
The Fourier expansion for $E_{k-2,2}(z,s+1-k,\xi)$ can be obtained in a convenient form using the Whittaker 
functions by applying the following proposition:
\begin{proposition}[Proposition 2.2, \cite{panchishkin:2003}]
Suppose that $s,l\in\mathbb{Z}$ are two integers satisfying $s\leqslant 0$ and $s+l>0$, 
then there is the following Fourier expansion:
\begin{equation*}
\begin{split}
&\mathbf{E}_{l,N}(z,s;a,b)=\delta\left(\frac{a}{N}\right)\,[\zeta(l+2s;b,N)+(-1)^{l+2s}\,\zeta(l+2s;-b,N)]\\
&\quad +\frac{(-2\pi i)^{l+2s}\,(-1)^s\,\Gamma(l+2s-1)}{(4\pi y)^{l+2s-1}\,N\,\Gamma(l+s)\,\Gamma(s)}\\
&\qquad \times [\zeta(l+2s-1;a,N)+(-1)^{l+2s}\,\zeta(l+2s-1;-a,N)]\\
&\quad +\frac{(-2\pi i)^{l+2s}\,(-1)^s}{N^{l+2s}\,\Gamma(l+s)}\\
&\qquad \times\left(\sum_{\substack{dd'>0\\d'\equiv a\,\mathrm{mod}\,N}}\mathrm{sgn}(d)\,d^{l+2s-1}\,e^{\frac{2\pi idb}{N}}\,W(\frac{4\pi dd'y}{N},l+s,s)\,e^{\frac{2\pi idd'z}{N}}\right)\,,
\end{split}
\end{equation*}
where $\delta(x)=1$ if $x$ is an integer, $\delta(x)=0$ otherwise, and
\begin{equation*}
\zeta(s;a,N)=\sum_{0<n\equiv a(\mathrm{mod}\,N)}n^{-s}
\end{equation*}
denote the partial Riemann zeta function (defined if necessary by analytic continuation on $s$).
\end{proposition}
In this Proposition the Eisenstein series $\mathbf{E}_{l,N}$ of weight $l$ are defined as
\begin{equation}
\mathbf{E}_{l,N}(z,s;a,b)=\sum(cz+d)^{-l}|cd+d|^{-2s}\,,
\end{equation}
where $(c,d)\equiv(a,b)\pmod{N}$ and $(c,d)\neq(0,0)$, $a,b\in\mathbb{Z}/N\mathbb{Z}$.

In order to apply the above proposition for our case, we note that
\begin{align}
&E_{18,2}(z,s+19,\xi)=\mathbf{E}_{18,2}(z,s-19,0,1)\,,\\
&k=20\,,\quad 1<s\leqslant 19\,,\quad a=0\,,\quad b=1\,,\quad N=2\,,\quad l=18\,,\\
&\delta(\frac{a}{N})=1\,,\\
&\zeta(s;0,2)=\sum_{0<n\equiv0(2)}n^{-s}=\sum_{n=1}^\infty(2n)^{-s}=2^{-s}\zeta(s)\,,\\
&\zeta(s;1,2)=\sum_{0<n\equiv1(2)}n^{-s}=\zeta(s)-\zeta(s;0,2)=(1-2^{-s})\zeta(s)\,.
\end{align}

The Whittaker function $W(y,\alpha,\beta)$ is defined as
\begin{equation}
W(y,\alpha,\beta)=\Gamma(\beta)^{-1}\int_0^{+\infty}(u+1)^{\alpha-1}\,u^{\beta-1}\,\mathrm{e}^{-yu}\,du
\end{equation}
for $y>0$, $\alpha,\beta\in\mathbb{C}$ with $\Re(\beta)>0$ and for arbitrary $\alpha$ and $\beta$ this 
function is defined by the analytic continuation and the functional equation:
\begin{equation}
W(y,\alpha,\beta)=y^{1-\alpha-\beta}\,W(y,1-\beta,1-\alpha)\,.
\end{equation}
For a non negative integer $r$, we have
\begin{equation}
W(y,\alpha,-r)=\sum_{i=0}^{r}\frac{(-1)^i\binom{r}{i}\Gamma(\alpha)}{\Gamma(\alpha-i)}\,y^{r-i}\,.
\end{equation}

Now we write the explicit expansion for the even weight $l=k-2$ (when $k$ is even) at the point $s+1-k$:
\begin{equation}
\label{eq:E-k-2}
\begin{split}
&(4\pi y)^{s+1-k}\,E_{k-2,2}(z,s+1-k,\xi)\\
&=(4\pi y)^{s+1-k}\,\Big(2\,\zeta(k-2+2(s+1-k);1,2)\\
&\quad+\frac{(-2\pi i)^{k-2+2(s+1-k)}\,(-1)^{s+1-k}\Gamma(k-2+2(s+1-k)-1)}{(4\pi y)^{k-2+2(s+1-k)-1}2\,\Gamma(k-2+s+1-k)\,\Gamma(s+1-k)}\\
&\qquad\times2\,\zeta(k-2+2(s+1-k)-1;0,2)\\
&\quad+\frac{(-2\pi i)^{k-2+2(s+1-k)}\,(-1)^{s+1-k}}{2^{k-2+2(s+1-k)}\Gamma(k-2+s+1-k)}\\
&\qquad\times\sum_{\substack{dd'>0\\d'\equiv 0\,\mathrm{mod}\,2}}\mathrm{sgn}(d)\,d^{k-2+2(s+1-k)-1}\,e^{\pi id}\\
&\qquad\times W(\frac{4\pi dd'y}{2},k-2+s+1-k,s+1-k)\,e^{\pi idd'z}\Big)\,.
\end{split}
\end{equation}
The first term of \eqref{eq:E-k-2} becomes
\begin{equation}
\begin{split}
&(4\pi y)^{s+1-k}\,2\,\zeta(k-2+2(s+1-k);1,2)\\
&\quad=2\,(4\pi y)^{s+1-k}\,(1-2^{k-2s})\,\zeta(2s-k)\,.
\end{split}
\end{equation}
The second term of \eqref{eq:E-k-2} becomes
\begin{equation}
\begin{split}
&(4\pi y)^{s+1-k}\,\frac{(2\pi i)^{2s-k}\,(-1)^{s+1-k}\Gamma(2s-k-1)}{(4\pi y)^{2s-k-1}2\,\Gamma(s-1)\,\Gamma(s+1-k)}\,2\,\zeta(2s-k-1;0,2)\\
&\quad=(4\pi y)^{2-s}\,\frac{(2\pi)^{2s-k}\,(-1)^{s-k/2}\,(-1)^{s+1-k}\,\Gamma(2s-k-1)}{\Gamma(s-1)\,\Gamma(s+1-k)}\times\\
&\qquad\times2^{k+1-2s}\,\zeta(2s-k-1)\\
&\quad=(4\pi y)^{2-s}\,\frac{2\,\pi^{2s-k}\,(-1)^{1-k/2}\,\Gamma(2s-k-1)\,\zeta(2s-k-1)}{\Gamma(s-1)\,\Gamma(s+1-k)}
\end{split}
\end{equation}
The last term of \eqref{eq:E-k-2} contains the sum over all $d$ and $d'$ such that the product $dd'$ is positive 
and $d'\equiv0\pmod{2}$.  Let us write $d'=2d_1$ and break down the summation in two parts, one for positive $d$ 
and another for negative $d$.  We also substitute $dd_1$ by $n$.  Therefore,
\begin{equation}
\begin{split}
&\sum_{\substack{dd'>0\\d'\equiv 0\,\mathrm{mod}\,2}}\mathrm{sgn}(d)\,d^{2s-k-1}\,e^{\pi id}\,W(2\pi dd'y,s-1,s+1-k)\,e^{\pi idd'z}\\
&=     \sum_{\substack{d>0\\d_1>0}}d^{2s-k-1}\,(-1)^d\,W(4\pi dd_1y,s-1,s+1-k)\,e^{2\pi idd_1z}\\
&\quad+\sum_{\substack{d<0\\d_1<0}}(-1)\,d^{2s-k-1}\,(-1)^d\,W(4\pi dd_1y,s-1,s+1-k)\,e^{2\pi idd_1z}\\
&=     \sum_{n>0}\sum_{d|n}d^{2s-k-1}\,(-1)^d\,W(4\pi ny,s-1,s+1-k)\,q^n\\
&\quad+\sum_{\substack{d>0\\d_1>0}}(-1)\,(-1)^{2s-k-1}\,d^{2s-k-1}\,(-1)^d\,W(4\pi dd_1y,s-1,s+1-k)\,e^{2\pi idd_1z}\\
&=     \sum_{n>0}\sum_{d|n}d^{2s-k-1}\,(-1)^d\,W(4\pi ny,s-1,s+1-k)\,q^n\\
&\quad+\sum_{n>0}\sum_{d|n}d^{2s-k-1}\,(-1)^d\,W(4\pi ny,s-1,s+1-k)\,q^n\\
&=2\sum_{n>0}\sum_{d|n}d^{2s-k-1}\,(-1)^d\,W(4\pi ny,s-1,s+1-k)\,q^n\,.
\end{split}
\end{equation}
The last term of \eqref{eq:E-k-2} becomes
\begin{equation}
\begin{split}
\frac{(-1)^{1-k/2}\,2\,\pi^{2s-k}}{\Gamma(s-1)}\sum_{n>0}\sum_{d|n}(-1)^d\,d^{2s-k-1}\,W(4\pi ny,s-1,s+1-k)\,q^n\,.
\end{split}
\end{equation}
To abbreviate the following manipulations we introduce some notations. Let
\begin{equation}
\label{eq:C-coeff}
\begin{split}
&C_0^\prime(s)=\frac{(-1)^{1-k/2}\,2\,\pi^{2s-k}\,\Gamma(2s-k-1)\,\zeta(2s-k-1)}{\Gamma(s-1)\,\Gamma(s+1-k)}\,,\\
&C_0^{\prime\prime}(s)=2\,(1-2^{k-2s})\,\zeta(2s-k)\,,\\
&C_1(s)=(-1)^{k/2}\,2\,\pi^{2s-k}\,,\\
&C_2(s)=(-1)^{k/2}\,(2-2^{2s-k})\,\pi^{2s-k}\,,\\
&C_3(s)=(-1)^{k/2}\,2\,(1+3^{2s-k-1})\,\pi^{2s-k}\,,\\
&C_4(s)=(-1)^{k/2}\,(2-2^{2s-k}-2^{4s-2k-1})\,\pi^{2s-k}\,,
\end{split}
\end{equation}
then \eqref{eq:E-k-2} can be written as
\begin{equation}
\begin{split}
&(4\pi y)^{s+1-k}\,E_{k-2,2}(z,s+1-k,\xi)\\
&\quad=C_0^\prime\,(4\pi y)^{2-s}+C_0^{\prime\prime}\,(4\pi y)^{s+1-k}\\
&\qquad+C_1\,\frac{W( 4\pi y,s-1,s+1-k)}{\Gamma(s-1)}(4\pi y)^{s+1-k}\,q\\
&\qquad+C_2\,\frac{W( 8\pi y,s-1,s+1-k)}{\Gamma(s-1)}(4\pi y)^{s+1-k}\,q^2\\
&\qquad+C_3\,\frac{W(12\pi y,s-1,s+1-k)}{\Gamma(s-1)}(4\pi y)^{s+1-k}\,q^3\\
&\qquad+C_4\,\frac{W(16\pi y,s-1,s+1-k)}{\Gamma(s-1)}(4\pi y)^{s+1-k}\,q^4+\cdots\,.
\end{split}
\end{equation}
To find the image of the projection operator as a $q$-expansion $\mathcal{H}ol(F)=\sum A_n(s)q^n$ we 
apply the Holomorphic Projection Lemma \cite[Proposition (5.1)]{gross-zagier:1986} 
(the Lemma is originally due to Sturm \cite{sturm:1980}). We denote the Fourier coefficients
of the product inside of the holomorphic projection operator by $\widetilde{A}_n(s,y)$:
\begin{align}
F(z,s,y) = G_{2,2}(z)(4\pi y)^{s+1-k}E_{k-2,2}(z,s+1-k,\xi )=\sum\widetilde{A}_n(s,y)q^n.
\end{align}
It should be noted, that the relevant polynomial decay hypotheses of the Lemma is satisfied for all actions on 
$E_{k-2,2}(z,s+1-k,\xi)\vert\gamma$ of $\gamma\in\mathrm{SL}_2(\mathbb{Z})$ and at each critical point $s$, 
see \cite[(2.3)]{panchishkin:2003}.

Consider the case of weight $k=20$.  Recall that the holomorphic projection belongs to the space of cups forms of 
weight $20$ and level $2$, which is a $4$-dimensional space. This space can be spanned by the the following four modular 
forms: $g_{20}(z)$, $g_{20}(2z)$ and the orthogonal subspace of newforms $h_{20,2}^{(1)}(z)$ and $h_{20,2}^{(2)}(z)$:
\begin{align}
&g_{20}(z) =q + 456\,q^2 + 50652\,q^3 - 316352\,q^4 - 2377410\,q^5 + \dotsc\,,\\
&g_{20}(2z)=q^2 + 456\,q^4 + 50652\,q^6 - 316352\,q^8 + \dotsc\,,\\
&h_{20,2}^{(1)}(z) = q - 512\,q^2 - 13092\,q^3 + 262144\,q^4 + 6546750\,q^5 + \dotsc\,,\\
&h_{20,2}^{(2)}(z) = q + 512\,q^2 - 53028\,q^3 + 262144\,q^4 - 5556930\,q^5 + \dotsc\,.
\end{align}
Consider the linear combination
\begin{equation}
\mathcal{H}ol(F)(z)=K_1\,g_{20}(z)+K_2\,g_{20}(2z)+K_3\,h_{20,2}^{(1)}(z)+K_4\,h_{20,2}^{(2)}(z)\,.
\end{equation}
Coefficients $K_i=K_i(s)$ can be found explicitly by comparing the Fourier coefficients (we also do not indicate 
the dependence of $A_i$ from $s$ to abbreviate the notation):
\begin{equation}
\label{eq:lin-syst-A1-A4}
\left\{
\begin{aligned}
A_1 &= K_1\cdot1         + K_2\cdot0   + K_3\cdot1        + K_4\cdot1\\
A_2 &= K_1\cdot456       + K_2\cdot1   + K_3\cdot(-512)   + K_4\cdot512\\
A_3 &= K_1\cdot50652     + K_2\cdot0   + K_3\cdot(-13092) + K_4\cdot(-53028)\\
A_4 &= K_1\cdot(-316352) + K_2\cdot456 + K_3\cdot262144   + K_4\cdot262144\,.
\end{aligned}
\right.
\end{equation}
Resolving the system of linear equations \eqref{eq:lin-syst-A1-A4}, we obtain
\begin{equation}
\label{eq:K}
\left\{
\begin{aligned}
K_{1} & = \frac{{-13 A_{4}} + {152 A_{3}} + {5928 A_{2}} + {8432992 A_{1}}}{22947840},\\
K_{2} & = \frac{{4229 A_{4}} + {24104 A_{3}} + {940056 A_{2}} - {311728736 A_{1}}}{2868480},\\
K_{3} & = \frac{{3 A_{4}} + {16 A_{3}} - {1368 A_{2}} + {762432 A_{1}}}{2039808},\\
K_{4} & = \frac{-A_{4} - {16 A_{3}} + {456 A_{2}} + {286144 A_{1}}}{1105920}\,.
\end{aligned}
\right.
\end{equation}
Recall that the computation of the Fourier coefficients of the Holomorphic projection is given by the formula 
\cite[Proposition (5.1)]{gross-zagier:1986}:
\begin{equation}
\label{eq:A-coeffs}
\begin{split}
A_m(s)&=\frac{(4\,\pi\,m)^{k-1}}{(k-2)!}\int_0^{\infty}\widetilde{A}_m(s,y)\,\mathrm{e}^{-4\pi my}\,y^{k-2}\,dy\\
&=\frac{1}{(k-2)!}\int_0^{\infty}\widetilde{A}_m(s,y)\,\mathrm{e}^{-4\pi my}\,(4\pi my)^{k-2}\,d(4\pi my)\,.
\end{split}
\end{equation}
The explicit computations of $A_m$ and $K_m$ are given in Appendix \ref{app:coeff}.

\section{Result for $L(s+10, g_{20})L(s+9, g_{20})$}
\label{sec:res2}

Recall that in each critical point $s\in \{-8,-6,-4,-2,0,1,3,5,7,9\}$ of the standard $L$-function 
$L(s,F_{12}, \mathrm{St})$ according to (\ref{eq3.1}) we compute two factors, one of which is 
$L(s+10, g_{20})L(s+9, g_{20})$. We consider the product $L(s, g_{20})L(s-1, g_{20})$  at points 
$s\in \{2,4,6,8,10,11,13,15,17,19\}$. In the Section \ref{sec:L(s,g20)L(s-1,g20)} we established the formula
\begin{align}
L(s, g_{20})L(s-1, g_{20})=\dfrac{L(s,g_{20}\otimes G_{2,2})}{1-456\cdot 2^{1-s}+2^{21-2s}}.
\end{align}
Later, we found that
\begin{align}
L(s,g_{20}\otimes G_{2,2})=\dfrac{3}{2}\dfrac{(4\pi)^{19}}{\Gamma(s)}\langle g_{20}(z), \mathcal{H}ol(F(z,s,y))\rangle,
\end{align}
where
\begin{align}
F(z,s,y)=G_{2,2}(z)(4\pi y)^{s-19}E_{18,2}(z,s-19,\xi),
\end{align}
We computed coefficients $K_1(s),K_2(s),K_3(s)$ and $K_4(s)$ so that
\begin{align}
\mathcal{H}ol(F)(z)=K_1g_{20}(z)+K_2g_{20}(2z)+K_3h^{(1)}_{20,2}(z)+K_4h^{(2)}_{20,2}(z).
\end{align}
Both forms $h^{(1)}_{20,2}(z)$ and $h^{(2)}_{20,2}(z)$ are orthogonal to $g_{20}$, therefore their Petersson inner 
product with $g_{20}$ is zero. We compute the Petersson inner product $\langle g_{20}(z),g_{20}(2z)\rangle$ using 
the trace operator (see \cite[paragraph 3.1]{serre:1973}).
\begin{align}
\begin{split}
\langle g_{20}(z),g_{20}(2z)\rangle & {}= 2^{-k/2}\cdot 3^{-1}\cdot 2^{1-k/2}\langle g_{20}(z),T_{2}(g_{20}(2z))\rangle\\
\\
&=2^{-16}\cdot 19\langle g_{20}(z),g_{20}(2z)\rangle.
\end{split}
\end{align}
Therefore,
\begin{align}
\langle g_{20}(z), \mathcal{H}ol(F(z,s,y))\rangle=(K_1(s)+2^{-16}\cdot 19\cdot K_2(s))\langle g_{20}(z),g_{20}(z)\rangle
\end{align}
The explicit values for each critical point $s$ of the product $L(s+10, g_{20})\,L(s+9, g_{20})$ are given in 
Appendix \ref{app:L(s,g20)L(s-1,g20)}.

\section{The main identity}

Combining the above results into the original expression (\ref{eq1.1}) we get the final result
\begin{align}
L(&s, F_{12}, \mathrm{St})=L(s+11,\Delta\otimes\Delta)\,L(s+10, g_{20})\,L(s+9, g_{20})\\
\nonumber\\
&=R_s^{(\Delta)}\,\pi_s^{(\Delta)}\,R_s^{(g_{20})}\,\pi_s^{(g_{20})}\,\langle\Delta,\Delta\rangle\,\langle g_{20}, g_{20}\rangle\\
\nonumber\\
&=R_s^{(\Delta)}\,\pi_s^{(\Delta)}\,\frac{3}{2}\,\frac{(4\pi)^{19}}{\Gamma(s+10)}\,\frac{(K_1(s)+2^{-16}\cdot 19\cdot K_2(s))}{1-456\cdot 2^{-9-s}+2^{1-2s}}\,\langle\Delta,\Delta\rangle\,\langle g_{20}, g_{20}\rangle.
\end{align}
where $R_s^{(\Delta)}, \pi_s^{(\Delta)}$ were computed in Section \ref{sec:Delta} 
and $K_1(s), K_2(s)$ were computed in Section \ref{sec:res2}. 
For each critical values $s\in\{-8, -6, -4, -2, 0, 1, 3, 5, 7, 9\}$ we evaluate 
the final expression in the form 
$L(s, F_{12}, \mathrm{St}) = R_s\,\pi_s\,\langle\Delta,\Delta\rangle\,\langle g_{20}, g_{20}\rangle$ 
with the rational coefficient $R_s$ and the corresponding power of $\pi$.
It might be interesting to see the values factorised in primes.

\medskip
{
\renewcommand{\arraystretch}{0.0}
\begin{tabular}{|r|cc|r|}
\hline&&\\\strut $s$ & $R_s$ & $\pi_s$ & numerical value\\\hline\rule{0pt}{6pt}&&&\\
-8 & $\dfrac{-1 \cdot 2^{31} \cdot 17 \cdot 11411 \cdot 1207259}{3 \cdot 5^2 \cdot 7 \cdot 11 \cdot 13 \cdot 61}$ & $\pi^{6}$ & -903525.807173
\\\rule{0pt}{6pt}&&&\\\hline\rule{0pt}{6pt}&&&\\
-6 & $\dfrac{-1 \cdot 2^{26} \cdot 47 \cdot 791797}{3^6 \cdot 17 \cdot 113}$ & $\pi^{12}$ & -14105.832863
\\\rule{0pt}{6pt}&&&\\\hline\rule{0pt}{6pt}&&&\\
-4 & $\dfrac{2^{24} \cdot 392033}{3^5 \cdot 5^3 \cdot 7 \cdot 17 \cdot 19}$ & $\pi^{18}$ & 728.260808
\\\rule{0pt}{6pt}&&&\\\hline\rule{0pt}{6pt}&&&\\
-2 & $\dfrac{-1 \cdot 2^{26} \cdot 479903}{3^8 \cdot 5^3 \cdot 7^3 \cdot 13 \cdot 17 \cdot 157}$ & $\pi^{24}$ & -24.122802
\\\rule{0pt}{6pt}&&&\\\hline\rule{0pt}{6pt}&&&\\
0  & $\dfrac{2^{22} \cdot 5779}{3^{13} \cdot 5^4 \cdot 7^3 \cdot 11 \cdot 13}$ & $\pi^{30}$ & 3.485667
\\\rule{0pt}{6pt}&&&\\\hline\rule{0pt}{6pt}&&&\\
1  & $\dfrac{2^{25} \cdot 2269}{3^{14} \cdot 5^4 \cdot 7^3 \cdot 11^2 \cdot 13 \cdot 17}$ & $\pi^{34}$ & 1.901053
\\\rule{0pt}{6pt}&&&\\\hline\rule{0pt}{6pt}&&&\\
3  & $\dfrac{2^{40}}{3^{16} \cdot 5^6 \cdot 7^4 \cdot 11^3 \cdot 13^2 \cdot 17}$ & $\pi^{42}$ & 1.156624
\\\rule{0pt}{6pt}&&&\\\hline\rule{0pt}{6pt}&&&\\
5  & $\dfrac{2^{40}}{3^{20} \cdot 5^8 \cdot 7^6 \cdot 11 \cdot 13^3 \cdot 17}$ & $\pi^{50}$ & 1.029466
\\\rule{0pt}{6pt}&&&\\\hline\rule{0pt}{6pt}&&&\\
7  & $\dfrac{2^{40}}{3^{23} \cdot 5^{10} \cdot 7^6 \cdot 11^2 \cdot 13^2 \cdot 17^2}$ & $\pi^{58}$ & 1.006025
\\\rule{0pt}{6pt}&&&\\\hline\rule{0pt}{6pt}&&&\\
9  & $\dfrac{2^{21} \cdot 9413 \cdot 6782351}{3^{23} \cdot 5^{10} \cdot 7^8 \cdot 11^4 \cdot 13^4 \cdot 17^2 \cdot 19 \cdot 61}$ & $\pi^{66}$ & 1.000909
\\\rule{0pt}{6pt}&&&\\\hline
\end{tabular}
}
{

\section{Numerical computation of the Petersson inner product}

To compute numerically the Petersson inner product of $\Delta$ by itself and $g_{20}$ 
by itself we use the classical result by Rankin (see \cite{rankin:1952}, Theorem 5):
\begin{align}
\langle f_k, f_k\rangle=\frac{(4\pi)^{1-k}(k-2)!}{\zeta(l)}\frac{\alpha_r}{\alpha_l+\alpha_r-\alpha_k}L(k-1, f_k)L(l, f_k),
\end{align}
where $f_k$ is the cusp form of weight $k=\{12, 14, 16, 18, 20\}$ of the form
\begin{align}
f_k(z)=E_{k-12}(z)\Delta(z)
\end{align}
and $4\leq r\leq k/2-2$;\\
$E_k$ denotes Eisenstein series
\begin{align*}
&E_k(z)=\sum_{n=0}^{\infty}\alpha_k(n)q^n,\\
&\alpha_k(0)=1,\\
&\alpha_k=\alpha_k(1)=-\frac{2k}{B_k},
\end{align*}
$B_k$ is a Bernoulli number.\\
For $f_k=\Delta$ we are able to use only one choice of critical value $l=8$. 
To compute the numerical values $L(11, \Delta)$ and $L(8, \Delta)$ we use 
Dokchitser's $L$-function Calculator. The obtained value is
\begin{align}
\langle\Delta,\Delta\rangle=0.000001035362056205680432094820996804,
\end{align}
which coincides with the values given by Zagier up to $11$ digit.\\
We used again Rankin's theorem to compute the Petersson inner product of $g_{20}$ by itself. 
For the modular form of weight $20$ there are three choices of $l=12, 14, 16$. 
For each choice of $l$ we computed the special value of $L(l, g_{20})$ using 
Dokchitser's $L$-function Calculator. The obtained values are
\begin{align}
\langle g_{20}, g_{20}\rangle=0.000008265541531659702744699575969 \mbox{ for } l=12,\\
\langle g_{20}, g_{20}\rangle=0.000008265541531659703390644766954,\\
\langle g_{20}, g_{20}\rangle=0.000008265541531659703069998511729.
\end{align}

\bibliographystyle{alpha}
\bibliography{Journal2pre}

\newpage
\appendix

\section{Explicit values for $L(s+11, \Delta\otimes\Delta)$}
\label{app:L(s,Delta2)}

The following are the values of $L(s+11,\Delta\otimes\Delta)$ presented in the form
$$L(s+11,\Delta\otimes\Delta)=R_s^{(\Delta)}\,\pi_s^{(\Delta)}\,\langle\Delta,\Delta\rangle$$

{
\renewcommand{\arraystretch}{0.0}
\begin{tabular}{|r|cc||r|cc||r|cc|}
\hline&&&&&&&&\\
\rule{0pt}{4pt}&&&&&&&&\\
$s$ & $R_s^{(\Delta)}$ & $\pi_s^{(\Delta)}$ & $s$ & $R_s^{(\Delta)}$ & $\pi_s^{(\Delta)}$ & $s$ & $R_s^{(\Delta)}$ & $\pi_s^{(\Delta)}$\\
\rule{0pt}{2pt}&&&&&&&&\\\hline\rule{0pt}{6pt}&&&&&&&&\\
-8 & $\dfrac{2^{20}\cdot9}{35}$ & $\pi^3$ & 0 & $\dfrac{2^{14}}{14175}$ & $\pi^{11}$ & 7 & $\dfrac{2^{35}}{15\cdot18!}$ &  $\pi^{25}$ \\
\rule{0pt}{6pt}&&&&&&&&\\\hline\rule{0pt}{6pt}&&&&&&&&\\
-6 & $\dfrac{-2^{16}}{9}$  & $\pi^5$ & 1 & $\dfrac{2^{23}}{11!}$ & $\pi^{13}$ & 9 & $\dfrac{2^{41}}{245\cdot20!}$ & $\pi^{29}$ \\
\rule{0pt}{6pt}&&&&&&&&\\\hline\rule{0pt}{6pt}&&&&&&&&\\
-4 & $\dfrac{2^{13}}{45}$ & $\pi^7$ & 3 & $\dfrac{2^{28}}{14!}$ & $\pi^{17}$ &&&\\
\rule{0pt}{6pt}&&&&&&&&\\\hline\rule{0pt}{6pt}&&&&&&&&\\
-2 & $\dfrac{-2^{14}}{2205}$ & $\pi^9$ & 5 & $\dfrac{2^{31}}{3\cdot16!}$ & $\pi^{21}$ &&&\\
\rule{0pt}{6pt}&&&&&&&&\\\hline
\end{tabular}
}


\newpage
\section{Explicit values of $C$, $A$ and $K$ coefficients}
\label{app:coeff}

First we compute the integrals given by \eqref{eq:A-coeffs}.

{\small
\begin{align}
\begin{split}
A_1{}&=\dfrac{C'_{0}}{(k-2)!}\int_{0}^{+\infty}(4\pi y)^{2-s}e^{-4\pi y}(4\pi y)^{k-2}d(4\pi y)\\
&+\dfrac{C''_{0}}{(k-2)!}\int_{0}^{+\infty}(4\pi y)^{s+1-k}e^{-4\pi y}(4\pi y)^{k-2}d(4\pi y)\\
&+\dfrac{C_1}{24(k-2)!}\int_{0}^{+\infty}\dfrac{W(4\pi y,s-1,s+1-k)}{\Gamma(s-1)}(4\pi y)^{s+1-k}e^{-4\pi y}(4\pi y)^{k-2}d(4\pi y)\\
&=\dfrac{\Gamma(k+1-s)}{(k-2)!}C'_0+\dfrac{\Gamma(s)}{(k-2)!}C''_0\\
&\quad{}+\dfrac{C_1}{24(k-2)!}\int_{0}^{+\infty}\sum_{i=0}^{k-s-1}\dfrac{(-1)^{i}\binom{k-s-1}{i}(4\pi y)^{k-2-i}}{\Gamma(s-1-i)}e^{-4\pi y}d(4\pi y)\\
&=\dfrac{\Gamma(k+1-s)}{(k-2)!}C'_0+\dfrac{\Gamma(s)}{(k-2)!}C''_0\\
&\quad{}+\dfrac{C_1}{24(k-2)!}\sum_{i=0}^{k-1-s}\dfrac{(-1)^{i}\binom{k-1-s}{i}\Gamma(k-1-i)}{\Gamma(s-1-i)}\\
\end{split}
\end{align}
}

{\small
\begin{align}
\begin{split}
A_2{}&=\dfrac{C'_{0}}{(k-2)!}\int_{0}^{+\infty}(4\pi y)^{2-s}e^{-8\pi y}(8\pi y)^{k-2}d(8\pi y)\\
&+\dfrac{C''_{0}}{(k-2)!}\int_{0}^{+\infty}(4\pi y)^{s+1-k}e^{-8\pi y}(8\pi y)^{k-2}d(8\pi y)\\
&+\dfrac{C_1}{(k-2)!}\int_{0}^{+\infty}\dfrac{W(4\pi y,s-1,s+1-k)}{\Gamma(s-1)}(4\pi y)^{s+1-k}e^{-8\pi y}(8\pi y)^{k-2}d(8\pi y)\\
&+\dfrac{C_2}{24(k-2)!}\int_{0}^{+\infty}\dfrac{W(8\pi y,s-1,s+1-k)}{\Gamma(s-1)}(4\pi y)^{s+1-k}e^{-8\pi y}(8\pi y)^{k-2}d(8\pi y)\\
&=\dfrac{C'_{0}}{(k-2)!}2^{s-2}\int_{0}^{+\infty}(8\pi y)^{k-s}e^{-8\pi y}d(8\pi y)\\
&\quad{}+\dfrac{C''_{0}}{(k-2)!}\int_{0}^{+\infty}(8\pi y)^{s+1-k}e^{-8\pi y}d(8\pi y)\\
&\quad{}+\dfrac{C_1}{(k-2)!}\int_{0}^{+\infty} \sum_{i=0}^{k-1-s}  \dfrac{(-1)^{i}\binom{k-1-s}{i}(4\pi y)^{k-1-s-i}}{\Gamma(s-1-i)}(4\pi y)^{s+1-k}e^{-8\pi y}(8\pi y)^{k-2}d(8\pi y) \\
&\quad{}+\dfrac{C_2}{24(k-2)!}\int_{0}^{+\infty} \sum_{i=0}^{k-1-s}  \dfrac{(-1)^{i}\binom{k-1-s}{i}(8\pi y)^{k-1-s-i}}{\Gamma(s-1-i)}(4\pi y)^{s+1-k}e^{-8\pi y}(8\pi y)^{k-2}d(8\pi y) \\
&=\dfrac{\Gamma(k+1-s)}{(k-2)!}2^{s-2}C'_{0}+\dfrac{\Gamma(s)}{(k-2)!}2^{k-1-s}C''_{0}\\
&\quad{}+\dfrac{C_1}{(k-2)!}\sum_{i=0}^{k-1-s} 2^{i}\dfrac{(-1)^{i}\binom{k-1-s}{i}\Gamma(k-1-i)}{\Gamma(s-1-i)}\\
&\quad{}+\dfrac{C_2}{24(k-2)!} 2^{k-1-s}\sum_{i=0}^{k-1-s} \dfrac{(-1)^{i}\binom{k-1-s}{i}\Gamma(k-1-i)}{\Gamma(s-1-i)}.
\end{split}
\end{align}
}

{\small
\begin{align}
\begin{split}
A_3{}&=\dfrac{4C'_{0}}{(k-2)!}\int_{0}^{+\infty}(4\pi y)^{2-s}e^{-12\pi y}(12\pi y)^{k-2}d(12\pi y)\\
&+\dfrac{4C''_{0}}{(k-2)!}\int_{0}^{+\infty}(4\pi y)^{s+1-k}e^{-12\pi y}(12\pi y)^{k-2}d(12\pi y)\\
&+\dfrac{C_1}{(k-2)!}\int_{0}^{+\infty}\dfrac{W(4\pi y,s-1,s+1-k)}{\Gamma(s-1)}(4\pi y)^{s+1-k}e^{-12\pi y}(12\pi y)^{k-2}d(12\pi y)\\
&+\dfrac{C_2}{(k-2)!}\int_{0}^{+\infty}\dfrac{W(8\pi y,s-1,s+1-k)}{\Gamma(s-1)}(4\pi y)^{s+1-k}e^{-12\pi y}(12\pi y)^{k-2}d(12\pi y)\\
&+\dfrac{C_3}{24(k-2)!}\int_{0}^{+\infty}\dfrac{W(12\pi y,s-1,s+1-k)}{\Gamma(s-1)}(4\pi y)^{s+1-k}e^{-12\pi y}(12\pi y)^{k-2}d(12\pi y)\\
&=\dfrac{4\cdot 3^{k-1-s}C'_{0}}{(k-2)!}\int_{0}^{+\infty}(12\pi y)^{k-s}e^{-12\pi y}d(12\pi y)\\
&\quad{}+\dfrac{4\cdot 3^{k-1-s}C''_{0}}{(k-2)!}\int_{0}^{+\infty}(12\pi y)^{s-1}e^{-12\pi y}d(12\pi y)\\
&\quad{}+\dfrac{C_1}{(k-2)!}\int_{0}^{+\infty} \sum_{i=0}^{k-1-s}  \dfrac{(-1)^{i}\binom{k-1-s}{i}(4\pi y)^{k-1-s-i}}{\Gamma(s-1-i)}(4\pi y)^{s+1-k}e^{-12\pi y}(12\pi y)^{k-2}d(12\pi y) \\
&\quad{}+\dfrac{C_2}{(k-2)!}\int_{0}^{+\infty} \sum_{i=0}^{k-1-s}  \dfrac{(-1)^{i}\binom{k-1-s}{i}(8\pi y)^{k-1-s-i}}{\Gamma(s-1-i)}(4\pi y)^{s+1-k}e^{-12\pi y}(12\pi y)^{k-2}d(12\pi y) \\
&\quad{}+\dfrac{C_3}{24(k-2)!}\int_{0}^{+\infty} \sum_{i=0}^{k-1-s}  \dfrac{(-1)^{i}\binom{k-1-s}{i}(12\pi y)^{k-1-s-i}}{\Gamma(s-1-i)}(4\pi y)^{s+1-k}e^{-12\pi y}(12\pi y)^{k-2}d(12\pi y) \\
&=\dfrac{\Gamma(k+1-s)}{(k-2)!}4\cdot 3^{s-2}C'_{0}+\dfrac{\Gamma(s)}{(k-2)!}4\cdot 3^{k-1-s}C''_{0}\\
&\quad{}+\dfrac{C_1}{(k-2)!}\sum_{i=0}^{k-1-s} 3^{i}\dfrac{(-1)^{i}\binom{k-1-s}{i}\Gamma(k-1-i)}{\Gamma(s-1-i)}\\
&\quad{}+\dfrac{2^{k-1-s}C_2}{(k-2)!} \sum_{i=0}^{k-1-s}\Big(\dfrac{3}{2}\Big)^{i} \dfrac{(-1)^{i}\binom{k-1-s}{i}\Gamma(k-1-i)}{\Gamma(s-1-i)}\\
&\quad{}+\dfrac{3^{k-1-s}C_3}{24(k-2)!} \sum_{i=0}^{k-1-s} \dfrac{(-1)^{i}\binom{k-1-s}{i}\Gamma(k-1-i)}{\Gamma(s-1-i)}.
\end{split}
\end{align}
}

{\small
\begin{align}
\begin{split}
A_4{}&=\dfrac{C'_{0}}{(k-2)!}\int_{0}^{+\infty}(4\pi y)^{2-s}e^{-16\pi y}(16\pi y)^{k-2}d(16\pi y)\\
&+\dfrac{C''_{0}}{(k-2)!}\int_{0}^{+\infty}(4\pi y)^{s+1-k}e^{-16\pi y}(16\pi y)^{k-2}d(16\pi y)\\
&+\dfrac{4C_1}{(k-2)!}\int_{0}^{+\infty}\dfrac{W(4\pi y,s-1,s+1-k)}{\Gamma(s-1)}(4\pi y)^{s+1-k}e^{-16\pi y}(16\pi y)^{k-2}d(16\pi y)\\
&+\dfrac{C_2}{(k-2)!}\int_{0}^{+\infty}\dfrac{W(8\pi y,s-1,s+1-k)}{\Gamma(s-1)}(4\pi y)^{s+1-k}e^{-16\pi y}(16\pi y)^{k-2}d(16\pi y)\\
&+\dfrac{C_3}{(k-2)!}\int_{0}^{+\infty}\dfrac{W(12\pi y,s-1,s+1-k)}{\Gamma(s-1)}(4\pi y)^{s+1-k}e^{-16\pi y}(16\pi y)^{k-2}d(16\pi y)\\
&+\dfrac{C_4}{24(k-2)!}\int_{0}^{+\infty}\dfrac{W(16\pi y,s-1,s+1-k)}{\Gamma(s-1)}(4\pi y)^{s+1-k}e^{-16\pi y}(16\pi y)^{k-2}d(16\pi y)\\
&=\dfrac{4^{s-2}C'_{0}}{(k-2)!}\int_{0}^{+\infty}(16\pi y)^{k-s}e^{-16\pi y}d(16\pi y)\\
&\quad{}+\dfrac{4^{k-1-s}C''_{0}}{(k-2)!}\int_{0}^{+\infty}(16\pi y)^{s-1}e^{-16\pi y}d(16\pi y)\\
&\quad{}+\dfrac{4C_1}{(k-2)!}\int_{0}^{+\infty} \sum_{i=0}^{k-1-s}  \dfrac{(-1)^{i}\binom{k-1-s}{i}(4\pi y)^{k-1-s-i}}{\Gamma(s-1-i)}(4\pi y)^{s+1-k}e^{-16\pi y}(16\pi y)^{k-2}d(16\pi y) \\
&\quad{}+\dfrac{C_2}{(k-2)!}\int_{0}^{+\infty} \sum_{i=0}^{k-1-s}  \dfrac{(-1)^{i}\binom{k-1-s}{i}(8\pi y)^{k-1-s-i}}{\Gamma(s-1-i)}(4\pi y)^{s+1-k}e^{-16\pi y}(16\pi y)^{k-2}d(16\pi y) \\
&\quad{}+\dfrac{C_3}{(k-2)!}\int_{0}^{+\infty} \sum_{i=0}^{k-1-s}  \dfrac{(-1)^{i}\binom{k-1-s}{i}(12\pi y)^{k-1-s-i}}{\Gamma(s-1-i)}(4\pi y)^{s+1-k}e^{-16\pi y}(16\pi y)^{k-2}d(16\pi y) \\
&\quad{}+\dfrac{C_4}{24(k-2)!}\int_{0}^{+\infty} \sum_{i=0}^{k-1-s}  \dfrac{(-1)^{i}\binom{k-1-s}{i}(16\pi y)^{k-1-s-i}}{\Gamma(s-1-i)}(4\pi y)^{s+1-k}e^{-16\pi y}(16\pi y)^{k-2}d(16\pi y) \\
&=\dfrac{\Gamma(k+1-s)}{(k-2)!}2^{2s-4}C'_{0}+\dfrac{\Gamma(s)}{(k-2)!}2^{2(k-1-s)}C''_{0}\\
&\quad{}+\dfrac{4C_1}{(k-2)!}\sum_{i=0}^{k-1-s} 2^{2i}\dfrac{(-1)^{i}\binom{k-1-s}{i}\Gamma(k-1-i)}{\Gamma(s-1-i)}\\
&\quad{}+\dfrac{2^{k-1-s}C_2}{(k-2)!} \sum_{i=0}^{k-1-s}2^{i} \dfrac{(-1)^{i}\binom{k-1-s}{i}\Gamma(k-1-i)}{\Gamma(s-1-i)}\\
&\quad{}+\dfrac{3^{k-1-s}C_3}{(k-2)!} \sum_{i=0}^{k-1-s}\Big(\dfrac{4}{3}\Big)^{i} \dfrac{(-1)^{i}\binom{k-1-s}{i}\Gamma(k-1-i)}{\Gamma(s-1-i)}\\
&\quad{}+\dfrac{2^{2(k-1-s)}C_4}{24(k-2)!} \sum_{i=0}^{k-1-s} \dfrac{(-1)^{i}\binom{k-1-s}{i}\Gamma(k-1-i)}{\Gamma(s-1-i)}.
\end{split}
\end{align}
}

Next, we give explicit values for coefficients $C$ from \eqref{eq:C-coeff} with 
$\pi$ power being factored out for brevity. All these coefficients share the same power 
of $\pi$, namely $\pi^{2(s+10)-k}=\pi^{2s}$ (the coefficient $C''_0$ has an implicit 
factor of $\pi$ through $\zeta$ function). Note, we need to take into account that the 
formulas were deduced for $L(s,g_{20})L(s-1,g_{20})$ where the $\pi$ factor has the form 
$\pi^{2s-k}$ while we need the values for $L(s+10,g_{20})L(s+9,g_{20})$, i.e. $s \mapsto s+10$.

\medskip
\noindent
\centerline
{\small
\renewcommand{\arraystretch}{0.0}
\begin{tabular}{|r|ccccccc|}
\hline&&&&&&&\\
\rule{0pt}{4pt}&&&&&&&\\
$s$ & $C'_0$				& $C''_0$ 				& $C_1$ & $C_2$ 			& $C_3$ 				& $C_4$ 				& $\pi^*$	\\\rule{0pt}{2pt}&&&&&&&\\\hline\rule{0pt}{6pt}&&&&&&&\\
-8  & $\dfrac{43867}{7182}$		& 0 					& 2 	& $\dfrac{131071}{65536}$ 	& $\dfrac{258280328}{129140163}$ 	& $\dfrac{17179738111}{8589934592}$ 	& $\pi^{-16}$	\\\rule{0pt}{6pt}&&&&&&&\\\hline\rule{0pt}{6pt}&&&&&&&\\
-6  & $\dfrac{35}{2}$			& 0 					& 2 	& $\dfrac{8191}{4096}$ 		& $\dfrac{3188648}{1594323}$ 		& $\dfrac{67100671}{33554432}$ 		& $\pi^{-12}$	\\\rule{0pt}{6pt}&&&&&&&\\\hline\rule{0pt}{6pt}&&&&&&&\\
-4  & $\dfrac{65}{6}$			& 0 					& 2 	& $\dfrac{511}{256}$ 		& $\dfrac{39368}{19683}$ 		& $\dfrac{261631}{131072}$ 		& $\pi^{-8}$	\\\rule{0pt}{6pt}&&&&&&&\\\hline\rule{0pt}{6pt}&&&&&&&\\
-2  & $\dfrac{11}{3}$			& 0 					& 2 	& $\dfrac{31}{16}$ 		& $\dfrac{488}{243}$ 			& $\dfrac{991}{512}$ 			& $\pi^{-4}$	\\\rule{0pt}{6pt}&&&&&&&\\\hline\rule{0pt}{6pt}&&&&&&&\\
0   & $\dfrac{3}{2}$ 			& 0 					& 2 	& 1 				& $\dfrac{8}{3}$ 			& $\dfrac{1}{2}$ 			& $\pi^0$	\\\rule{0pt}{6pt}&&&&&&&\\\hline\rule{0pt}{6pt}&&&&&&&\\
1   & 0 				& $\dfrac{1}{4}$ 			& 2 	& -2 				& 8 					& -10 					& $\pi^2$	\\\rule{0pt}{6pt}&&&&&&&\\\hline\rule{0pt}{6pt}&&&&&&&\\
3   & 0 				& $\dfrac{1}{480}$ 			& 2 	& -62 				& 488 					& -2110 				& $\pi^6$	\\\rule{0pt}{6pt}&&&&&&&\\\hline\rule{0pt}{6pt}&&&&&&&\\
5   & 0 				& $\dfrac{31}{1451520}$ 		& 2 	& -1022 			& 39368 				& -525310 				& $\pi^{10}$	\\\rule{0pt}{6pt}&&&&&&&\\\hline\rule{0pt}{6pt}&&&&&&&\\
7   & 0 				& $\dfrac{5461}{24908083200}$		& 2 	& -16382 			& 3188648 				& -134234110 				& $\pi^{14}$	\\\rule{0pt}{6pt}&&&&&&&\\\hline\rule{0pt}{6pt}&&&&&&&\\
9   & 0 				& $\dfrac{3202291}{1422749712384000}$	& 2 	& -262142 			& 258280328 				& -34360000510 				& $\pi^{18}$	\\\rule{0pt}{6pt}&&&&&&&\\\hline
\end{tabular}
}
\bigskip

\newpage
Substituting these values into the formulas for coefficients $A$ we obtain their explicit values at each point $s$.
The common $\pi$ factor is the same as for coefficients $C$.

\medskip
\noindent
\centerline
{\small
\renewcommand{\arraystretch}{0.0}
\begin{tabular}{|r|cccc|}
\hline&&&&\\
\rule{0pt}{4pt}&&&&\\
$s$ & $A_1$					& $A_2$ 				& $A_3$ 				& $A_4$ 					\\\rule{0pt}{2pt}&&&&\\\hline\rule{0pt}{6pt}&&&&\\
-8  & $\dfrac{88931}{14364}$ 			& $\dfrac{157008449}{14364}$ 		& $\dfrac{39586640915}{3591}$ 		& $\dfrac{24277850760593}{14364}$ 		\\\rule{0pt}{6pt}&&&&\\\hline\rule{0pt}{6pt}&&&&\\
-6  & $\dfrac{71}{1224}$ 			& $\dfrac{8387}{306}$ 			& $\dfrac{447871}{34}$ 			& $\dfrac{131485894}{153}$ 			\\\rule{0pt}{6pt}&&&&\\\hline\rule{0pt}{6pt}&&&&\\
-4  & $\dfrac{1}{6528}$ 			& $\dfrac{173}{3672}$ 			& $\dfrac{5103}{544}$ 			& $\dfrac{9380}{17}$ 				\\\rule{0pt}{6pt}&&&&\\\hline\rule{0pt}{6pt}&&&&\\
-2  & $\dfrac{1}{3144960}$ 			& $\dfrac{1}{8568}$ 			& $\dfrac{23801}{4455360}$ 		& $\dfrac{11015}{41769}$ 			\\\rule{0pt}{6pt}&&&&\\\hline\rule{0pt}{6pt}&&&&\\
0   & $\dfrac{1}{784143360}$ 			& $\dfrac{23}{49008960}$ 		& $\dfrac{4997}{196035840}$ 		& $\dfrac{421}{2042040}$ 			\\\rule{0pt}{6pt}&&&&\\\hline\rule{0pt}{6pt}&&&&\\
1   & $\dfrac{1}{5292967680}$ 			& $\dfrac{53}{1323241920}$ 		& $\dfrac{27}{5445440}$ 		& $\dfrac{-23}{33081048}$ 			\\\rule{0pt}{6pt}&&&&\\\hline\rule{0pt}{6pt}&&&&\\
3   & $\dfrac{199}{1270312243200}$ 		& $\dfrac{49}{5671036800}$ 		& $\dfrac{2059}{3207859200}$ 		& $\dfrac{-10529}{4962157200}$ 			\\\rule{0pt}{6pt}&&&&\\\hline\rule{0pt}{6pt}&&&&\\
5   & $\dfrac{19}{65330343936}$ 		& $\dfrac{1277}{285820254720}$ 		& $\dfrac{8167}{63515612160}$ 		& $\dfrac{-6631}{10508097600}$ 			\\\rule{0pt}{6pt}&&&&\\\hline\rule{0pt}{6pt}&&&&\\
7   & $\dfrac{286703}{400148356608000}$ 	& $\dfrac{28267}{10003708915200}$ 	& $\dfrac{633}{15247232000}$ 		& $\dfrac{-8745697}{6252318072000}$ 		\\\rule{0pt}{6pt}&&&&\\\hline\rule{0pt}{6pt}&&&&\\
9   & $\dfrac{4803437}{2134124568576000}$ 	& $\dfrac{4737913}{2134124568576000}$ 	& $\dfrac{20552747}{533531142144000}$ 	& $\dfrac{-7037087527}{2134124568576000}$ 	\\\rule{0pt}{6pt}&&&&\\\hline
\end{tabular}
}
\bigskip
%

\newpage
Finally, we list the explicit values for the coefficients $K$ computed by expression \eqref{eq:K}.
The common $\pi$ factor is the same as for coefficients $C$ and $A$ above.
Since the coefficients $K_3$ and $K_4$ do not contribute to the final expression, we do ommit them here.

\medskip
\noindent
\centerline
{\small
\renewcommand{\arraystretch}{0.0}
\begin{tabular}{|r|cc|}
\hline&&\\
\rule{0pt}{4pt}&&\\
$s$ & $K_1$						& $K_2$ 						\\\rule{0pt}{2pt}&&\\\hline\rule{0pt}{6pt}&&\\
-8  & $\dfrac{-435883731901}{495673344}$ 		& $\dfrac{3045934023523439}{1177224192}$ 		\\\rule{0pt}{6pt}&&\\\hline\rule{0pt}{6pt}&&\\
-6  & $\dfrac{-217211831}{585169920}$ 			& $\dfrac{100968174943}{73146240}$ 			\\\rule{0pt}{6pt}&&\\\hline\rule{0pt}{6pt}&&\\
-4  & $\dfrac{-255571}{1404407808}$ 			& $\dfrac{156430715}{175550976}$ 			\\\rule{0pt}{6pt}&&\\\hline\rule{0pt}{6pt}&&\\
-2  & $\dfrac{45173}{1369297612800}$ 			& $\dfrac{74862131}{171162201600}$ 			\\\rule{0pt}{6pt}&&\\\hline\rule{0pt}{6pt}&&\\
0   & $\dfrac{36097}{56232488632320}$ 			& $\dfrac{3748999}{7029061079040}$ 			\\\rule{0pt}{6pt}&&\\\hline\rule{0pt}{6pt}&&\\
1   & $\dfrac{23831}{210871832371200}$ 		& $\dfrac{876017}{26358979046400}$ 			\\\rule{0pt}{6pt}&&\\\hline\rule{0pt}{6pt}&&\\
3   & $\dfrac{4553}{69773768064000}$ 			& $\dfrac{-1256}{105304870125}$ 			\\\rule{0pt}{6pt}&&\\\hline\rule{0pt}{6pt}&&\\
5   & $\dfrac{424061}{3881958732288000}$ 		& $\dfrac{-1672}{55749637125}$ 				\\\rule{0pt}{6pt}&&\\\hline\rule{0pt}{6pt}&&\\
7   & $\dfrac{923549}{3483809118720000}$ 		& $\dfrac{-66896}{850539335625}$ 			\\\rule{0pt}{6pt}&&\\\hline\rule{0pt}{6pt}&&\\
9   & $\dfrac{8127882069959}{9794709827950215168000}$ 	& $\dfrac{-304138734083887}{1224338728493776896000}$ 	\\\rule{0pt}{6pt}&&\\\hline
\end{tabular}
}

\newpage
\section{Explicit values of the product\\ $L(s+10, g_{20})\,L(s+9, g_{20})$}
\label{app:L(s,g20)L(s-1,g20)}

{
\renewcommand{\arraystretch}{0.0}
\begin{tabular}{|r|cc|}
\hline&&\\
\rule{0pt}{4pt}&&\\
$s$	& $R_s^{(g_{20})}$					& $\pi_s^{(g_{20})}$\\\rule{0pt}{2pt}&&\\\hline\rule{0pt}{6pt}&&\\
-8	& $\dfrac{-479626345744384}{1177605}$			& $\pi^3$	\\\rule{0pt}{6pt}&&\\\hline\rule{0pt}{6pt}&&\\
-6	& $\dfrac{38107606016}{155601}$				& $\pi^7$	\\\rule{0pt}{6pt}&&\\\hline\rule{0pt}{6pt}&&\\
-4	& $\dfrac{802883584}{1526175}$				& $\pi^{11}$	\\\rule{0pt}{6pt}&&\\\hline\rule{0pt}{6pt}&&\\
-2	& $\dfrac{1965682688}{4426469775}$			& $\pi^{15}$	\\\rule{0pt}{6pt}&&\\\hline\rule{0pt}{6pt}&&\\
0	& $\dfrac{1479424}{3447969525}$				& $\pi^{19}$	\\\rule{0pt}{6pt}&&\\\hline\rule{0pt}{6pt}&&\\
1	& $\dfrac{2323456}{175846445775}$			& $\pi^{21}$	\\\rule{0pt}{6pt}&&\\\hline\rule{0pt}{6pt}&&\\
3	& $\dfrac{8388608}{145073317764375}$			& $\pi^{25}$	\\\rule{0pt}{6pt}&&\\\hline\rule{0pt}{6pt}&&\\
5	& $\dfrac{16777216}{34367988873684375}$		& $\pi^{29}$	\\\rule{0pt}{6pt}&&\\\hline\rule{0pt}{6pt}&&\\
7	& $\dfrac{2097152}{436209089550609375}$		& $\pi^{33}$	\\\rule{0pt}{6pt}&&\\\hline\rule{0pt}{6pt}&&\\
9	& $\dfrac{63842269963}{1305893808013068186412500}$	& $\pi^{37}$	\\\rule{0pt}{6pt}&&\\\hline
\end{tabular}
}

\end{document}